\documentclass[11pt]{article}
\usepackage{amssymb}
\usepackage{amsfonts}
\usepackage{amsthm}
\usepackage{graphicx}
\usepackage[utf8]{inputenc}
\usepackage{xcolor}
\usepackage[english]{babel}
\usepackage{amsmath}
\usepackage{eucal}
\usepackage{enumerate}
\usepackage{stmaryrd}
\usepackage{amsfonts,epsfig,epstopdf,titling,url,array}
\usepackage{verbatim} 
\usepackage{hyperref}

\numberwithin{equation}{section}

\newtheorem{thm}{Theorem}
\newtheorem{prop}{Proposition}
\newtheorem{lemma}{Lemma}

\newcommand{\beq}{\begin{equation}}
\newcommand{\RN}[1]{%
  \textup{\uppercase\expandafter{\romannumeral#1}}%
}
\newcommand{\eeq}{\end{equation}}

\newcommand{\beqs}{\begin{equation*}}
\newcommand{\eeqs}{\end{equation*}}

\usepackage[
top    = 2.75cm,
bottom = 2.50cm,
left   = 2.49cm,
right  = 2.549cm]{geometry}

\title{Stability of the solutions in a predator-prey model with density-dependent diffusion 
}

\author{
Leoncio Rodriguez Q., Luis F. Gordillo\\ Department of Mathematics and Statistics \\ Utah State University \\
Logan, UT 84322.
}

\begin{document}
\bibliographystyle{plain}
\large{

\begin{titlepage}
\LARGE{

\vskip 1truein

\begin{center}
Stability of steady state solutions in a predator-prey model with density-dependent diffusion 
\end{center}
}
\large{
\vskip 1truein
\begin{center} 
Leoncio Rodriguez Q.\footnote{ Corresponding author, e-mail address: leoncio.quinones@usu.edu. ORCID-ID: 0000-0003-1615-1973}, Luis F. Gordillo.\footnote{e-mail address:  luis.gordillo@usu.edu} \\
Department of Mathematics and Statistics\\
Utah State University \\
Logan, UT 84322\\
\end{center}
\date{}

\vskip 1truein
\begin{abstract}
In this note we present a study of the solutions associated to a particular spatial extension of the Rosenzweig-MacArthur model for predator and prey. The analysis presented here shows that positive steady state solutions emerge via a transcritical bifurcation mechanism, in accordance with the insight obtained from previous numerical and analytical results. In the model under discussion, prey is assumed to move avoiding crowds via a density-dependent diffusion and also incorporates the existence of a refuge zone, where predators cannot consume prey. Saturation in prey consumption is also included through a Holling type II functional response.
\end{abstract}
}
\end{titlepage}

\maketitle

\section{Introduction}
Rosenzweig and MacArthur's work \cite{RosenzweigMac} gave origin to the popular model bearing their names, which is currently included in virtually any textbook of theoretical biology. The model has been extensively studied and applied in multiple situations, see for instance \cite{Kot,Turchin}. A particular extension of the Rosenzweig-MacArthur model that includes spatial movement of predators and prey is studied in \cite{Leo1}, where the existence of non-trivial positive steady states is established. Specifically, the model in \cite{Leo1} assumes that (1) prey spread follows a nonlinear diffusion rule, (2) prey have a refuge zone (or “protection zone”) where predators cannot enter, and (3) predators move following linear diffusion. This scenario is relevant to understand how density-dependent diffusion in the prey might affect the dynamics, in contrast with a system that exhibits simple linear diffusion. In our case, the nonlinear diffusion might be interpreted as prey avoiding crowds, a possible mechanism to avoid predators.

The analysis elaborated in \cite{Leo1}, which is complemented with careful numerical simulations, involves an adaptation of a maximum principle to the nonlinear diffusion case and makes use of classical results in bifurcation theory for partial differential equations to show the existence of positive solutions. Those developments, however, lack the analysis of the bifurcation mechanism that is finally observed. In terms of stability, the interpretation of the results remains speculative, only associated to the underlying biological explanations. 

In this paper we provide an analytical proof that the bifurcation found in \cite{Leo1} corresponds to the transcritical case and show that the nontrivial solutions form one of the stable branches, thus closing the gap in the results previously obtained.

\section{The model}
Here we summarize the proposed extension of the Rosenzweig-MacArthur equations given in \cite{Leo1}, which are defined over a bounded domain $\Omega\subset\mathbb{R}^2$ that represents a closed environment where predators and preys live. We consider an additional domain, called the ``refuge zone", $\Omega_{0}\subset\Omega$, where  predators cannot enter. We assume that $\Omega$  and $\Omega_0$ have sufficiently smooth boundaries, that $\overline{\Omega}_{0}\subset\Omega$, and define $\Omega_{1}=\Omega\setminus\overline{\Omega}_{0}$. The system of parabolic equations for the prey and predator populations, denoted by $u$ and $v$ respectively, is given by 
\begin{equation}\label{system0}
\begin{aligned}
    \partial_{t}u &= D_u\nabla\cdot u\nabla u + ru\left(1 - \frac{u}{\lambda}\right) - \frac{b(x)uv}{1+mu}\qquad\text{in}\quad \Omega, \\ 
    \partial_{t}v &= D_v\Delta v - \mu v + \frac{cuv}{1+mu} \qquad\qquad \text{in}\quad \Omega_{1},\\
    v&\equiv 0\qquad \text{in}\quad \Omega\setminus\Omega_1,
\end{aligned}
\end{equation}
with boundary and initial conditions
\begin{equation}
\begin{aligned}
    \partial_{n}u &= 0 \quad\text{on}\quad \partial\Omega, \\
    \partial_{n}v &= 0 \quad\text{on}\quad \partial\Omega_{1},\\ 
    u(x,0) &= u_{0}(x)\geq 0 \quad\text{for}\quad x \in \Omega, \\
    v(x,0) &= v_{0}(x)\geq 0 \quad\text{for}\quad x \in \Omega_{1}. 
\end{aligned}
\end{equation}
The parameters are positive and the function $b(x)$, which determines the efficiency of predator attacks, is defined by
\begin{equation}
b(x) = 
    \begin{cases}
    b > 0& \qquad\text{if } x \in \Omega_{1},\\
    0& \qquad \text{if } x \in \overline{\Omega}_{0},
    \end{cases}
\end{equation}
thus characterizing the refuge zone $\Omega_{0}$. The non-flux boundary condition on $\partial\Omega_{1}$ restricts predators to the exterior of the refuge zone. However, prey can move freely over the whole domain $\Omega$.
As in \cite{Leo1}, we consider the dimensionless form of the model (2.1) at the steady state, given by
\begin{equation}\label{system1}
\begin{aligned}
    0 &= \nabla\cdot u\nabla u +\lambda u - u^{2} - \frac{b(x)uv}{1+mu}\qquad\text{in}\quad \Omega, \\ 
    0 &= \Delta v - \mu v + \frac{cuv}{1+mu}\qquad\text{in}\quad \Omega_{1}.\\ 
\end{aligned}
\end{equation}
It is shown in \cite{Leo1} that positive steady state solutions for this problem exist and the solution curves are contrasted numerically with the case where prey move with linear diffusion. 

\section{Preliminary results}

The main result in \cite{Leo1} (Theorem 1) establishes the existence of the curve of positive solutions for the system (2.4). This curve is determined by 
\begin{align}\label{dos}
    \{(\mu(s),u(s),v(s))=(\mu_{\lambda}(s),\lambda-s\alpha_{\mu_{\lambda}}(x)+o(|s|),s+o(|s|)):s\in(0,a)\},
\end{align}
which was obtained by using a well known result due to Crandall and Rabinowitz, see \cite{crandall}. This curve also satisfies the following conditions, $\mu_{\lambda}(0)=c\lambda/(1+m\lambda)$, $u(0)=\lambda$, $v(0)=0$ and $F(\mu_{\lambda}(0),\lambda,0)=0$, where $F$ is the operator defined by Equation (2.2) in \cite{Leo1}, that is

\begin{equation}\label{operator1}
F(\mu,w,v)=\begin{pmatrix}
-\nabla\cdot w\nabla w + \lambda\Delta w - \lambda w + w^{2} + \frac{b(x)(\lambda-w)v}{1+m(\lambda-w)}\\ \\
\Delta v - \mu v + \frac{c(\lambda - w)v}{1+m(\lambda-w)}
\end{pmatrix}^{T}.   
\end{equation}

In what follows we shall assume that $F:U\times V \rightarrow Z$ is $C^{2}(U\times V,Z)$, where $U\subseteq \mathbb{R}$, $V\subseteq X = X_{\Omega}\times X_{\Omega_{1}}$, $Z = Y_{\Omega}\times Y_{\Omega_{1}}$ 
and  
\begin{align*}
 X_{\Omega} &=\{u\in W^{2,p}(\Omega):\partial_{n}u=0 \text{ on } \partial\Omega\}, \  Y_{\Omega}=L^{p}(\Omega),\\
 X_{\Omega_{1}} &= \{u\in W^{2,p}(\Omega_{1}):\partial_{n}u=0 \text{ on } \partial\Omega_{1}\},\ Y_{\Omega_{1}}=L^{p}(\Omega_{1}),
\end{align*}
so  $(\mu_{\lambda}(0),u(0),v(0))\in U\times V \subseteq \mathbb{R}\times X$. By defining $w$ as $w=\lambda-u$, the expression $F_{(w,v)}(\mu_{\lambda}(0),0,0)$, determined by 

\begin{equation}\label{operator3}
    F_{(w,v)}(\mu,0,0)[\alpha,\beta]=\begin{pmatrix}
    \lambda\Delta\alpha - \lambda\alpha + \frac{b(x)\lambda\beta}{1+m\lambda}\\ \\
    \Delta\beta - \mu\beta + \frac{c\lambda\beta}{1+m\lambda}
    \end{pmatrix}^{T},
\end{equation}
is equivalent to $F_{(u,v)}(\mu_{\lambda}(0),\lambda,0)$ and therefore (see \cite{Leo1})
\[dim(N(F_{(w,v)}(\mu_{\lambda}(0),0,0)))=1 \qquad\textnormal{and}\qquad  dim(N(F_{(u,v)}(\mu_{\lambda}(0),\lambda,0)))=1,\]
where $N(T)$ denotes the kernel of a linear operator $T$. Consequently, we only need to consider the operator $F_{(u,v)}(\mu_{\lambda}(0),\lambda,0)$.

Let us denote 
\[R=Range(F_{(u,v)}(\mu_{\lambda}(0),\lambda,0)) \quad\text{and}\quad N=N(F_{(u,v)}(\mu_{\lambda}(0),\lambda,0)),\] 
and define the operator $G:U\times V\times (R\cap X)\times \mathbb{R}\rightarrow Z$ by
\begin{align}\label{tres}
    G(\mu,u,v,\widetilde{w},\gamma) = F_{(u,v)}(\mu,u,v)(\Tilde{v}_{0}+\widetilde{w})-\gamma(\Tilde{v}_{0}+\widetilde{w}),
\end{align}
where  $N(F_{(u,v)}(\mu_{\lambda}(0),\lambda,0))=span\{(\alpha_{\mu_{\lambda}},1)\}=span\{\Tilde{v}_{0}\}$ so $\Tilde{v}_{0}\not\in R$ and $\Tilde{v}_{0}=(\alpha_{\mu_{\lambda}},1)$. We notice  that the Banach space $X$ admits a decomposition $X = N\oplus (R\cap X)$ (see equation $I.7.6$, pag 22 in \cite{Hansjorg}) and that $R\subset Z$, $X\subset Z$, and $R\cap X\subset X\hookrightarrow Z$. We emphasize that $\widetilde{w}$ must belong to $R\cap X$ and notice $\Tilde{v}_{0}\not\in R$ since we know that $\int_{\Omega_{1}}1dx\neq 0$, see \cite{Leo1}. Also, because $F_{(u,v)}(\mu,u,v):X\rightarrow Z$ and $\Tilde{v}_{0}+\widetilde{w}\in N\oplus (R\cap X)$,  we have that (\ref{tres}) is well defined. Furthermore, $G(\mu_{\lambda}(0),\lambda,0,0,0)= F_{(u,v)}(\mu_{\lambda}(0),\lambda,0)[\Tilde{v}_{0}]=0$.

\begin{lemma}\label{lema1}
The derivatives of $G$, defined in (\ref{tres}), are 
\begin{equation}\label{parteuno}
\begin{aligned}
    G_{\widetilde{w}}(\mu,u,v,\widetilde{w},\gamma)&=F_{(u,v)}(\mu,u,v)-\gamma I,\\ G_{\gamma}(\mu,u,v,\widetilde{w},\gamma)&=-(\Tilde{v}_{0}+\widetilde{w}),\\
    and \ \ G_{(\widetilde{w},\gamma)}(\mu,u,v,\widetilde{w},\gamma)&=F_{(u,v)}(\mu,u,v)-\gamma I-(\Tilde{v}_{0}+\widetilde{w}),
\end{aligned}{} 
\end{equation}

where $I$ is the identity operator on $X$.
\end{lemma}
\begin{proof}
First, we compute the first variation of $G$ with respect to $\widetilde{w}$. Consider
\begin{align}\label{cuatro}
    G(\mu,u,v,\widetilde{w}+\epsilon\chi,\gamma)=F_{(u,v)}(\mu,u,v)(\Tilde{v}_{0}+\widetilde{w}+\epsilon\chi)-\gamma(\Tilde{v}_{0}+\widetilde{w}+\epsilon\chi),
\end{align}
with $\chi\in R\cap X$. Since
\begin{equation}\label{cinco}
\begin{aligned}
    \frac{d}{d\epsilon}\left(F_{(u,v)}(\mu,u,v)(\Tilde{v}_{0}+\widetilde{w}+\epsilon\chi)\right)\Big|_{\epsilon=0}&= \frac{d}{d\epsilon}\left(F_{(u,v)}(\mu,u,v)(\Tilde{v}_{0}+\widetilde{w})\right)+ \frac{d}{d\epsilon}\left(F_{(u,v)}(\mu,u,v)(\epsilon\chi)\right)\\
    &=F_{(u,v)}(\mu,u,v)(\chi),
\end{aligned}
\end{equation}
we get 
\begin{equation}\label{seis}
    \begin{aligned}
    G_{\widetilde{w}}(\mu,u,v,\widetilde{w},\gamma)(\chi)&= \frac{d}{d\epsilon}\left(F_{(u,v)}(\mu,u,v)(\Tilde{v}_{0}+\widetilde{w}+\epsilon\chi)\right)\Big|_{\epsilon=0} - \gamma\chi\\
    &=F_{(u,v)}(\mu,u,v)(\chi) - \gamma\chi,
\end{aligned}{}
\end{equation}
and the first derivative in (\ref{parteuno}) follows from (\ref{seis}). For the second expression in (\ref{parteuno}), consider
\begin{align}
    G(\mu,u,v,\widetilde{w},\gamma+\epsilon\theta) = F_{(u,v)}(\mu,u,v)(\Tilde{v}_{0}+\widetilde{w})-(\gamma + \epsilon\theta)(\Tilde{v}_{0}+\widetilde{w}). 
\end{align}
Then
\begin{align}\label{siete}
    G_{\gamma}(\mu,u,v,\widetilde{w},\gamma)(\theta)=\frac{d}{d\epsilon}\left(G(\mu,u,v,\widetilde{w},\gamma+\epsilon\theta)\right)\Big|_{\epsilon=0} = -\theta(\Tilde{v}_{0}+\widetilde{w})
\end{align}
and the second equation in  (\ref{parteuno})  follows from (\ref{siete}). The last equation in (\ref{parteuno}) is obtained in a similar way, by calculating
\begin{align}
\frac{d}{d\epsilon}\left(G(\mu,u,v,\widetilde{w}+\epsilon\chi,\gamma+\epsilon\theta)\right)\Big|_{\epsilon=0}.
\end{align}
\end{proof}

By the previous Lemma,
\begin{equation}
    \begin{aligned}
        G_{\widetilde{w}}(\mu_{\lambda}(0),\lambda,0,0,0)&=F_{(u,v)}(\mu_{\lambda}(0),\lambda,0),\\
        G_{\gamma}(\mu_{\lambda}(0),\lambda,0,0,0)&=-\Tilde{v}_{0},\\
        G_{(\widetilde{w},\gamma)}(\mu_{\lambda}(0),\lambda,0,0,0)&=F_{(u,v)}(\mu_{\lambda}(0),\lambda,0)-\Tilde{v}_{0}.
    \end{aligned}
\end{equation}
In order to determine the stability of the solutions for our problem we need the following result that guarantees the existence of a differentiable curve of perturbed eigenvalues for the linearized operator.

\begin{prop}\label{propo1}
Assume that the operator $F$ satisfies $F\in C^{2}(U\times V,Z)$, where $U,V$,$R$ and $Z$ are defined as above, and that $(\mu_{\lambda}(0),\lambda,0)\in U \times V \subseteq \mathbb{R}\times X$. Then, there exist a continuously differentiable curve of perturbed eigenvalues, \[\{\gamma(s):s\in(-\delta,\delta), \gamma(0)=0\}\] 
such that  
\begin{equation}\label{pertuevals}
F_{(u,v)}(\mu(s),u(s),v(s))(\Tilde{v}_{0}+\widetilde{w}(s))=\gamma(s)(\Tilde{v}_{0}+\widetilde{w}(s)).   
\end{equation}
In this sense, $\gamma(s)$ is the perturbation of the zero (simple) eigenvalue of $F_{(u,v)}(\mu_{\lambda}(0),\lambda,0)$.
\end{prop}{}

\begin{proof}
Let $(\chi,\theta)\in N(G_{(\widetilde{w},\gamma)}(\mu_{\lambda}(0),\lambda,0,0,0))$, i.e. $G_{(\widetilde{w},\gamma)}(\mu_{\lambda}(0),\lambda,0,0,0)[\chi,\theta]=0$, which implies $F_{(u,v)}(\mu_{\lambda}(0),\lambda,0)(\chi)=\theta\Tilde{v}_{0}$. Thus, since $\theta\in\mathbb{R}$, we get $\theta\Tilde{v}_{0}\in R$. Suppose that $\theta\neq0$, then $\Tilde{v}_{0}\in R$. But we showed that $N(F_{(u,v)}(\mu_{\lambda}(0),\lambda,0))=span\{\Tilde{v}_{0}\}$, thus implying $\Tilde{v}_{0}\not\in R$, which is a contradiction (see also $I.7.4$ page $21$ in \cite{Hansjorg}), therefore we must have $\theta=0$. Now, if $\theta=0$ then $F_{(u,v)}(\mu_{\lambda}(0),\lambda,0)(\chi)=0$ and, consequently, $\chi=a\Tilde{v}_{0}$ for some constant $a$. Notice that $(\chi,\theta)\in N\times\{0\}$ and so $N(G_{(\widetilde{w},\gamma)}(\mu_{\lambda}(0),\lambda,0,0,0))\subseteq N\times\{0\}$. If $\chi\in N$ then $\chi\not\in R$, so if $\chi\in N\cap R$ then we obtain $\chi = 0$.

If we consider the operator $G_{(\widetilde{w},\gamma)}(\mu_{\lambda}(0),\lambda,0,0,0):(R\cap X)\times\mathbb{R}\rightarrow Z$ then it is clear that \[N(G_{(\widetilde{w},\gamma)}(\mu_{\lambda}(0),\lambda,0,0,0))\subseteq \times(R \cap X)\times\mathbb{R},\]
and then 
\[N(G_{(\widetilde{w},\gamma)}(\mu_{\lambda}(0),\lambda,0,0,0))\subseteq (N\cap R)\times\{0\}=\{0\}_{X}\times\{0\},\] 
implying that $G_{(\widetilde{w},\gamma)}(\mu_{\lambda}(0),\lambda,0,0,0)$ is invertible. Now, consider $\Tilde{z}\in Z$, which can be written as $\Tilde{z}=x+y$, where $x\in N(G_{(\widetilde{w},\gamma)}(\mu_{\lambda}(0),\lambda,0,0,0))$ and $y\in R\cap X$ (see eq. $I.7.5$ page $22$, \cite{Hansjorg}). Then \[G_{(\widetilde{w},\gamma)}(\mu_{\lambda}(0),\lambda,0,0,0)[\xi,\theta]=F_{(u,v)}(\mu_{\lambda}(0),\lambda,0)(\xi)-\theta\Tilde{v}_{0} = x + y,\] 
so we choose $\xi$ and $\theta$ such that $F_{(u,v)}(\mu_{\lambda}(0),\lambda,0)(\xi)=y$ and $x=-\theta\Tilde{v}_{0}$. Then the operator $G_{(\widetilde{w},\gamma)}(\mu_{\lambda}(0),\lambda,0,0,0)$ is onto, and it is not hard to see that it is also linear. Hence, from the argument above, we have that \[G_{(\widetilde{w},\gamma)}(\mu_{\lambda}(0),\lambda,0,0,0):(R\cap X)\times\mathbb{R}\rightarrow Z\] is an isomorphism. By the Implicit Function Theorem, (see Thm $I.4.1$ Pag 12, \cite{Hansjorg}), there exist differentiable functions $\widetilde{w}:U_{1}\times V_{1}\rightarrow R\cap X$ and  $\gamma:U_{1}\times V_{1}\rightarrow \mathbb{R}$ such that $(\lambda,0)\in V_{1}\subset X$, $\mu_{\lambda}\in U_{1}\subset \mathbb{R}$, $\widetilde{w}(\mu_{\lambda}(0),\lambda,0)=0$, $\gamma(\mu_{\lambda}(0),\lambda,0)=0$ and $G(\mu,u,v,\widetilde{w}(\mu,u,v),\gamma(\mu,u,v))=0$ for all $(\mu,u,v)\in U_{1}\times V_{1}$. By inserting the curves from ($\ref{dos}$) into $\widetilde{w}$ and $\gamma$, we obtain 
\begin{equation}\label{ocho}
    \begin{aligned}{}
    \gamma(s) &= \gamma(\mu(s),u(s),v(s))=\gamma(\mu_{\lambda}(s),\lambda-s\alpha_{\mu_{\lambda}}(x)+o(|s|),s+o(|s|)),\\
    \widetilde{w}(s) &=\widetilde{w}(\mu(s),u(s),v(s))=\widetilde{w}(\mu_{\lambda}(s),\lambda-s\alpha_{\mu_{\lambda}}(x)+o(|s|),s+o(|s|)), \ \ s\in(0,\delta)
    \end{aligned}
\end{equation}{}
for some $\delta>0$. Notice that $\widetilde{w}(0)=0$ and $\gamma(0)=0$. From (\ref{ocho}) and the condition \[G(\mu,u,v,\widetilde{w}(\mu,u,v),\gamma(\mu,u,v))=0\] we get
\begin{equation}
F_{(u,v)}(\mu(s),u(s),v(s))(\Tilde{v}_{0}+\widetilde{w}(s))=\gamma(s)(\Tilde{v}_{0}+\widetilde{w}(s)).   \end{equation}
\end{proof}{}


\section{Main result}
We now focus on the question whether the spectrum of $F_{(u,v)}(\mu_{\lambda}(0),\lambda,0)$ is in the left complex plane. We look at the eigenvalue equation
\begin{equation}
F_{(u,v)}(\mu_{\lambda}(0),\lambda,0)[\alpha,\beta]=\Lambda[\alpha,\beta],    
\end{equation}
where the operator in the left-hand side is defined as
\begin{equation*}
    F_{(w,v)}(\mu,0,0)[\alpha,\beta]=\begin{pmatrix}
    \lambda\Delta\alpha - \lambda\alpha + \frac{b(x)\lambda\beta}{1+m\lambda}\\ \\
    \Delta\beta - \mu\beta + \frac{c\lambda\beta}{1+m\lambda}
    \end{pmatrix}^{T},
\end{equation*}
see equation (2.4) in \cite{Leo1}. Thus, we have the following system of PDEs with boundary conditions
\begin{equation}\label{nueve}
    \begin{aligned}{}
    \Delta\alpha-\alpha + \frac{b(x)\beta}{1+m\lambda}&=\Lambda\alpha \ \ x\in\Omega \\
    \Delta\beta -\mu\beta+ \frac{c\lambda\beta}{1+m\lambda}&=\Lambda\beta \ \ x\in\Omega_{1}\\
    \partial_{n}\alpha &= 0, \ \ x\in\partial\Omega,\\
    \partial_{n}\beta &= 0, \ \ x\in\partial\Omega_{1}.
    \end{aligned}
\end{equation}
From the second equation in (\ref{nueve}), we have the Neumann eigenvalue problem
\begin{equation}
    \begin{aligned}{}
    -\Delta\beta &= \left(-\Lambda-\mu+\frac{c\lambda}{1+m\lambda}\right)\beta \ \ x\in\Omega_{1}\\
    \partial_{n}\beta &= 0, \ \ x\in\partial\Omega_{1}.
    \end{aligned}
\end{equation}
In general, the Neumann eigenvalues of the negative Laplacian are non-negative and therefore $-\mu+c\lambda/(1+m\lambda)\geq\Lambda$. However, in order for $\beta$ to not change sign, we consider $\mu=-\Lambda+c\lambda/(1+m\lambda)$, which corresponds to the zero eigenvalue (and hence $\beta$ is the positive constant eingenfunction assoticiated to it). Since our interest is the case $\Lambda<0$, we must assume $\Lambda\leq-\mu+c\lambda/(1+m\lambda)<0$. By the principle of linearized stability, we must require that $c\lambda/(1+m\lambda)=\mu_{\lambda}(0)<\mu$.
To determine the stability of solutions of (\ref{system1}) we also need to determine the sign of the perturbed eigenvalues $\gamma(s)$ for small values of $s\in(0,\delta)$. Since $\mu^{\prime}_{\lambda}(0)<0$, the bifurcation is \textit{transcritical} (see page 18, \cite{Hansjorg}), i.e. we have two solution curves intersecting at the bifurcation point:  $(\mu_{\lambda},\lambda,0)$, namely the semitrivial solution line $\Gamma_{u}=\{(\mu,u,v)=(\mu,\lambda,0):\mu>0\}$ and  the curve 
of nontrivial solutions $\{(\mu(s),u(s),v(s))\}$ determined by
\begin{align}
\{(\mu,u,v)=(\mu_{\lambda}(s),\lambda-s\alpha_{\mu_{\lambda}}(x)+o(|s|),s+o(|s|)):s\in(0,a)\}
 \end{align}
and obtained in Theorem 1 in \cite{Leo1} using known results on bifurcation from simple eigenvalues and theory of elliptic PDE's  (see \cite{crandall},\cite{lintakagi},\cite{lopezgomez},\cite{louni},\cite{pucciserrin},\cite{rabinowitz},\cite{shiwang},\cite{He2017}). 

Consider equation (\ref{pertuevals}) and its parametrization near $\mu_{\lambda}(0)$ (we use the parameter $r$) so it becomes
\begin{equation}
    F_{(u,v)}(r,\lambda,0)(\Tilde{v}_{0}+\widetilde{w}(r))=\gamma(r)(\Tilde{v}_{0}+\widetilde{w}(r)).
\end{equation}{}
By taking the derivative respect to $r$ we get
\begin{equation}\label{diez}
    F_{r(u,v)}(r,\lambda,0)(\Tilde{v}_{0}+\widetilde{w}(r))+F_{(u,v)}(r,\lambda,0)\dot{\widetilde{w}}(r)=\dot{\gamma}(r)(\Tilde{v}_{0}+\widetilde{w}(r))+\gamma(r)\dot{\widetilde{w}}(r),
\end{equation}
where $\dot{}=d/dr$. Since we are considering the parameter $r$, we write $\gamma(0)=\gamma(\mu_{\lambda}(0),\lambda,0)=\gamma(\mu_{\lambda}(0))=0$ and $\widetilde{w}(0)=\widetilde{w}(\mu_{\lambda}(0),\lambda,0)=\widetilde{w}(\mu_{\lambda}(0))=0$. Then, at $r=\mu_{\lambda}(0)$,  (\ref{diez}) becomes
\begin{equation}\label{once}
F_{r(u,v)}(\mu_{\lambda}(0),\lambda,0)\Tilde{v}_{0} + F_{(u,v)}(\mu_{\lambda}(0),\lambda,0)\dot{\widetilde{w}}(\mu_{\lambda}(0))=\dot{\gamma}(\mu_{\lambda}(0))\Tilde{v}_{0}
\end{equation}
and we can choose an element $\Tilde{v}^{*}_{0}\in Z^{*}$ (the dual of $Z$) such that $\langle \Tilde{v}_{0},\Tilde{v}^{*}_{0}\rangle =1$ and $ \langle z,\Tilde{v}^{*}_{0}\rangle = 0$ for any $z\in R$, (see $I.7.8$, page $22$ in \cite{Hansjorg}). The brackets $\langle \cdot,\cdot \rangle$ are used to denote the duality pairing between $Z$ and its dual. By applying the duality pairing to (\ref{once}) with respect to the element $\Tilde{v}^{*}_{0}$ and noticing that since $ F_{(u,v)}(\mu_{\lambda}(0),\lambda,0)\dot{\widetilde{w}}(\mu_{\lambda}(0))\in R$, then $\langle F_{(u,v)}(\mu_{\lambda}(0),\lambda,0)\dot{\widetilde{w}}(\mu_{\lambda}(0)), \Tilde{v}^{*}_{0}\rangle = 0$. We obtain the following expression of the derivative of the perturbed simple eigenvalue at the parameter value $\mu_{\lambda}(0)$,
\begin{equation}\label{doce}
   \langle F_{r(u,v)}(\mu_{\lambda}(0),\lambda,0)\Tilde{v}_{0}, \Tilde{v}^{*}_{0}\rangle = \dot{\gamma}(\mu_{\lambda}(0)).
 \end{equation}
Notice that if $F_{r(u,v)}(\mu_{\lambda}(0),\lambda,0)\Tilde{v}_{0}\not\in R$ then $\dot{\gamma}(\mu_{\lambda}(0))\neq 0$, and if $\mu>\mu_{\lambda}(0)$ then 
\begin{equation}\label{trece}
    \dot{\gamma}(\mu_{\lambda}(0))= \langle F_{r(u,v)}(\mu_{\lambda}(0),\lambda,0)\Tilde{v}_{0}, \Tilde{v}^{*}_{0}\rangle = -\int_{\Omega_{1}}dx = |\Omega_{1}|<0
\end{equation}
determines loss of stability. Therefore, the semitrivial solution curve $\Gamma_{u}$ is locally stable for $\mu>\mu_{\lambda}(0)$ and unstable for $\mu<\mu_{\lambda}(0)$.
For the bifurcating solution curve $\{\mu(s),u(s),v(s)\}$ we know that $\mu^{\prime}_{\lambda}(0)<0$, (see equation (2.19) in \cite{Leo1}), and
\begin{equation}\label{catorce}
    \langle F_{(u,v)(u,v)}(\mu_{\lambda}(0),\lambda,0)[\Tilde{v}_{0},\Tilde{v}_{0}],\Tilde{v}^{*}_{0}\rangle +2\mu^{\prime}_{\lambda}(0)\dot{\gamma}(\mu_{\lambda}(0))=0.
\end{equation}
Also (see $I.7.22$ page 24 in \cite{Hansjorg}),
\begin{equation}\label{quince}
    \langle F_{(u,v)(u,v)}(\mu_{\lambda}(0),\lambda,0)[\Tilde{v}_{0},\Tilde{v}_{0}],\Tilde{v}^{*}_{0}\rangle +\mu^{\prime}_{\lambda}(0)\dot{\gamma}(\mu_{\lambda}(0))=\dot{\widehat{\gamma}}(0),
\end{equation}
where $\widehat{\gamma}$ is the eigenvalue perturbation given in (\ref{pertuevals}) along the curve (and not $\gamma(\mu_{\lambda}(s))$, which is the eingenvalue perturbation along the function $\mu_{\lambda}(s)$). By subtracting (\ref{quince}) from  (\ref{catorce}) we obtain
\begin{equation}\label{diezyseis}
    \mu^{\prime}_{\lambda}(0)\dot{\gamma}(\mu_{\lambda}(0))=-\dot{\widehat{\gamma}}(0).
\end{equation}
We have that $\dot{\gamma}(\mu_{\lambda}(0))<0$, which means that there is a loss of stability in the semi-trivial solution $\Gamma_{u}$ at $(\mu_{\lambda}(0),\lambda,0)$, where  $\mu(0)=\mu_{\lambda}(0)$, and $\widehat{\gamma}(0)=0$ (since at $0$ is the simple zero eigenvalue). We know that $\mu^{\prime}_{\lambda}(0)<0$ and $\dot{\widehat{\gamma}}(0)<0$ then $\widehat{\gamma}(s)<0$ for $s>0$. Therefore, if $\mu^{\prime}_{\lambda}(0)<0$ then the conditions $\mu(s)<\mu_{\lambda}(0)$ and $\widehat{\gamma}(s)<0$ must hold simultaneously for $s>0$. Similarly, $\mu(s)>\mu_{\lambda}(0)$ and $\widehat{\gamma}(s)>0$ must hold for $s<0$. In other words, 
\begin{equation}\label{diezysiete}
    sign(\mu(s)-\mu_{\lambda}(0))=sign(\widehat{\gamma}(s)), \ \ s\in(-\delta,\delta).
\end{equation}
Therefore, the curve $\{(\mu(s),u(s),v(s))\}$ is stable whenever $\mu(s)<\mu_{\lambda}(0)$ and unstable if $\mu(s)>\mu_{\lambda}(0)$. Summarizing, the following result on the stability of the steady state solutions for the system (\ref{system1}) holds, 
\begin{thm}
Consider the curve of nontrivial positive solutions \begin{align}\label{diezyocho}
    \mathcal{C}(s)=\{(\mu(s),u(s),v(s))=(\mu_{\lambda}(s),\lambda-s\alpha_{\mu_{\lambda}}(x)+o(|s|),s+o(|s|)):s\in(0,a)\},
\end{align}
and the curve of semi-trivial solutions $\Gamma_{u}$
for the system (\ref{system1}), 
satisfying the conditions $\mu_{\lambda}(0)=c\lambda/(1+m\lambda)$, $u(0)=\lambda$, $v(0)=0$, and $F(\mu_{\lambda}(0),\lambda,0)=0$, where $F:U\times V \rightarrow Z$ is $C^{2}(U\times V,Z)$ and $U\subseteq \mathbb{R}$. Then, $\mathcal{C}(s)$ and $\Gamma_{u}$ satisfy the properties
\begin{enumerate}
    \item $\mathcal{C}(s)$ is stable whenever $\mu(s)<\mu_{\lambda}(0)$ and unstable whenever $\mu(s)>\mu_{\lambda}(0)$.
    \item $\Gamma_{u}$ is stable whenever $\mu(s)>\mu_{\lambda}(0)$ and unstable whenever $\mu(s)<\mu_{\lambda}(0)$
\end{enumerate}
and therefore the parameter value $\mu_{\lambda}(0)$ determines a transcritical bifurcation.
\end{thm}{}

\section{Discussion}
In this paper we continue the analysis of a nonlinear diffusion mechanism introduced in a prey population and initiated in \cite{Leo1}. The model exhibits a plausible prey adaptation response that counteracts predation. In addition, the spatial domain for the model contains a refuge zone for the prey, which excludes the presence of predators, and predator saturation on prey consumption is included in the equations through a Holling type II functional response.

Previously in \cite{Leo1}, it was determined that positive solutions exist at the steady state, a conclusion relying directly on the application of Crandall and Rabinowitz bifurcation results. Although a biological interpretation could give an intuitive idea of the solutions' stability nature, an analytical argument was lacking. This paper closes that gap by presenting a detailed theoretical argument, which shows that the mechanism involved corresponds to that of a transcritical bifurcation. The results obtained by the consideration of nonlinear diffusion complement recent studies that include the simultaneous effects of nonlinearities and prey refuge, see for instance \cite{wang-li,He2017,zhangrongzhang}.

\bibliography{REFS}
}

\end{document}